\documentclass[12pt]{article}
\usepackage{amsfonts}
\usepackage{graphicx}
\usepackage{amsmath}
\usepackage{hyperref}
\usepackage[cp866]{inputenc}
\newtheorem{defn}{Definition}[section]
\newtheorem{thm}{Theorem}[section]

\newtheorem{rem}{Remark}[section]

\newtheorem{prop}{Proposition}[section]
\begin{document}

\begin{center}
\textbf{A FRACTIONAL GENERALIZATION
OF THE POISSON PROCESSES AND SOME OF ITS PROPERTIES}
\vskip.3cm
\textbf{Nicy Sebastian} \\
Indian Statistical Institute, Chennai Centre, Taramani, Chennai - 600113, India\\
nicy@isichennai.res.in\\
\vskip.2cm
\textbf{Rudolf Gorenflo}\\
Department of Mathematics and Informatics, Free University of Berlin, Arnimallee 3, D-14195 Berlin, Germany
\end{center}

\begin{abstract}
We have provided a fractional generalization of the Poisson renewal processes
by replacing the first time derivative in the relaxation equation of the
survival probability by a fractional derivative of order $\alpha ~(0 < \alpha  \leq 1)$. A generalized Laplacian model associated with the Mittag-Leffler distribution is examined. We also discuss some properties of this new model and its relevance to time series. Distribution of gliding sums, regression behaviors
and sample path properties are studied. Finally we introduce the $q$-Mittag-Leffler process associated with the $q$-Mittag-Leffler distribution.
\end{abstract}
\section{Introduction}
\label{intro}
It is our intention to provide via fractional calculus a generalization of
the pure and compound Poisson processes, which are known to play a
fundamental role in renewal theory.
If the waiting time is exponentially distributed we have a Poisson
process, which is Markovian. However, other waiting time distributions are
also relevant in applications, in particular such ones with a fat tail caused by
a power law decay of its density. In this context we analyze a non-Markovian
renewal process with a waiting time distribution described by the Mittag-Leffler function. This distribution, containing the exponential as particular
case, is known to play a fundamental role in the infinite thinning procedure
of a generic renewal process governed by a power-asymptotic waiting time.

The concept of renewal process has been developed as a stochastic model
for describing the class of counting processes for which the times between
successive events are independently and identically distributed non-negative
random variables, obeying a given probability law. These times are referred
to as waiting times or inter-arrival times.

For a renewal process having waiting times $T_1, T_2, . . .,$ let $$t_0= 0,~~~ t_k= \sum_{j=1}^{k}T_j,~~~ k\geq1,$$
where $t_1 = T_1$ is the time of the first renewal, $t_2 = T_1 + T_2$ is the time of
the second renewal and in general $t_k$ denotes the $k^{th}$ renewal. The process is specified if we know the probability law for the waiting
times.
A relevant
quantity is the counting function $N(t)$ defined as
$$ N(t) = \max \left\{k|t_k \leq t, k = 0, 1, 2,\ldots \right\} , $$
that represents the effective number of events before or at instant $t$.
Also
 $$F_k(t)= P(t_k= T_1+ T_2+\cdots+T_k \leq t), f_k(t)= \frac{d}{d t}F_k(t), k \geq 1,$$
where $F_k(t)$ represents the probability that the sum of the first $k$ waiting
times is less or equal $t$ and $f_k(t)$ its density. We assume the waiting times $T_j=t_j-t_{j-1}$ to be mutually independent, all having the same probability density $f(t)$. We introduce the function $R(t)= P(T>t)= \int_t^\infty f(t'){\rm d}t',$ the survival probability. This name comes from theory of  maintenance and means the probability that the relevant piece of equipment lives at least until instant $t$. We have $$f(t)=-\frac{d}{dt}R(t).$$ In the classical Poisson process the survival probability obeys the relaxation equation
\begin{equation}\label{1.1}
\frac{d}{d t}R(t)= -\lambda R(t), t\geq 0; R(0+)=1, \lambda>0
 \end{equation}with a positive constant $\lambda$. We get

$$
R(t)=\exp(-\lambda t), t\geq 0,~~~ and~~~~ f(t)=\lambda {\exp}(-\lambda t), \lambda > 0, t\geq 0.
$$
Without going into details see \cite{Ross :1996}.
We know the probability that $k$ events occur in the interval of length $t$
is given by the well-known Poisson distribution
$$
P(N(t)=k)=\frac{(\lambda t)^k}{k!}{\exp}(-\lambda t),  t\geq 0, k=0,1,2,\ldots.
$$
Its mean is given as, $E(N(t))= \lambda t$.
The probability density $f_k(t)$ for the sum $t_k= T_1+ T_2+\cdots+T_k $ is the $k$-fold convolution of the waiting time density. In our case we have
$$
f_k(t)=\lambda \frac{(\lambda t)^{k-1}}{(k-1)!}{\exp}(-\lambda t),  t\geq 0, k = 1,2,\ldots,$$
so that the Erlang distribution function of order $k$ turns out to be
$$
F_k(t)=\int_0^t f_k(t')dt'= 1-\sum_{n=0}^{k-1}\frac{(\lambda t)^{n}}{n!}{\exp}(-\lambda t)=\sum_{n=k}^{\infty}\frac{(\lambda t)^{n}}{n!}{\exp}(-\lambda t), t\geq 0.
$$

A ``fractional" generalization of the Poisson renewal process is simply
obtained by generalizing the differential equation (\ref{1.1}) replacing there the
first derivative with the integro-differential operator $_tD_{*}^{\alpha}$
that is interpreted
as the fractional derivative of order $\alpha$ in Caputo's sense (\cite{Podlubny:1999}, \cite{Gorenflo and Mainardi:1997}) and which in the case $0<\alpha\leq1$ is defined as follows:
$$\label{eq:66.1}
_tD_{*}^{\alpha}f(t)=\left\{\begin{array}{ll}
f'(t)
 & if ~\alpha=1\\
 \frac{1}{\Gamma(1-\alpha)}\int_0^t \frac{f'(\tau)}{(t-\tau)^{\alpha}} {\rm d}\tau &if~ 0<\alpha<1.
 \end{array} \right .
$$
We write, taking for simplicity $\lambda = 1$,
\begin{equation}\label{1.7}
_tD_{*}^{\alpha}R(t)= - R(t), t>0, 0<\alpha \leq 1; R(0+)=1.
 \end{equation}
The solution of (\ref{1.7}) is known to be,
$$
R(t)= E_\alpha(-t^\alpha), t\geq 0, 0<\alpha \leq 1,
$$
where $E_\alpha(\cdot)$ denotes the Mittag-Leffler function (see \cite{Erdelyi: Magnus:Oberhettinger:Tricomi:1955}) which is given as the case $\beta=1$ of the
 two-index Mittag-Leffler function is defined as $$E_{\alpha,\beta}(z)=\sum_{k=0}^\infty\frac{ z^k}
{\Gamma(\beta+\alpha k)}, z\in \mathcal{C}.$$
 In contrast to the Poissonian case $ \alpha= 1,$ in the case $0 < \alpha< 1$ for large $t$
the function $R(t)$ no longer decays exponentially but algebraically.
As a consequence of the power-law asymptotics the process turns to be no
longer Markovian but of long-memory type. However, we recognize that
for $0 < \alpha< 1$  the function 	$R(t)$, keeps the completely monotonic
character of the Poissonian case. Complete monotonicity of a functions
$g(t)$ means
$$(-1)^n \frac{d^n}{dt^n}g(t)\geq 0, ~~~n=0,1,2,\ldots, ~~~ t\geq0,$$ or equivalently,  its representability as real Laplace transform of nonnegative
generalized function (or measure). By using the Laplace transform technique we generalize the Poisson
distribution to the fractional Poisson distribution,
$$P(N(t)=k)=\frac{t^{k\alpha}}{k!}E_\alpha^{(k)}(-t^\alpha), k=0, 1, 2, \ldots.$$
The corresponding fractional Erlang pdf (of order $k \geq 1$) is
$$f_k(t)=\alpha\frac{t^{k\alpha-1}}{(k-1)!}E_\alpha^{(k)}(-t^\alpha),$$
and the fractional Erlang distribution function turns out to be
$$\int_0^t f_k(t')dt'= 1-\sum_{n=0}^{k-1}\frac{(t)^{n\alpha}}{n!}E_\alpha^{(n)}(-t^\alpha)=
\sum_{n=k}^{\infty}\frac{(t)^{n \alpha}}{n!}E_\alpha^{(n)}(-t^\alpha).$$

\begin{figure}\begin{center}
\resizebox*{6cm}{!}{ \includegraphics{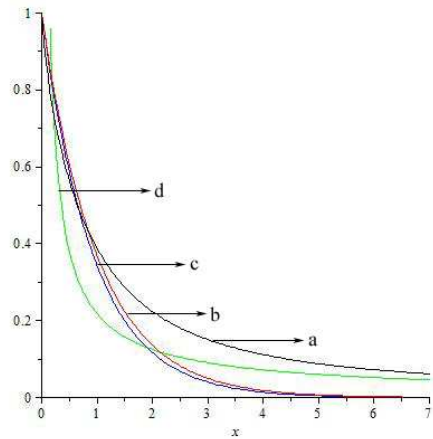}}\end{center}
\caption{(a)~
Mittag-Leffler decay for  $\alpha=0.8,$~
(b)~standard exponential decay  for
$\alpha=1,$~(c)~stretched exponential decay  for
$\alpha=0.8,$~(d)~power law decay  for
$\alpha=0.8$\label{fig:0}}
\end{figure}

When solving certain problems in processes of decay and oscillation, diffusion and wave propagation the solution can often be
obtained in terms of exponential or logarithmic
functions when the orders of differentiation are integer numbers. In the case of non-integer orders Mittag-Leffler functions
 or more general special functions are required, see  \cite{Gorenflo and Mainardi:2009}, \cite{Mathai:Saxena:Haubold:2006}.
 Pillai \cite{Pillai:1990} proved that $F_\alpha (x)=
1-E_\alpha(-x^\alpha),~ 0< \alpha \leq 1,~x>0$ and $F_\alpha (x)=0$
for $x\leq0$ are distribution functions, having the Laplace
transform $(1+s^{\alpha})^{-1}, s>0$.  He called $ F_\alpha (x)$, for
$0<\alpha\leq1,$ a Mittag-Leffler distribution and showed that $1-F_{\alpha}(x)$ is completely monotone for $x>0$.  For $x\in \mathcal{R}$ and $\alpha=1,$ the Mittag-Leffler function with argument $-x^{\alpha}$ reduces to a standard exponential decay $\exp(-x);$ when $0<\alpha<1,$ the Mittag-Leffler function is approximated for small values of $x$ by a stretched exponential decay (Weibull function) ${\exp}(-x^{\alpha}/\Gamma(\alpha+1))$ and for large values of $x$ by a power law $bx^{-\alpha}$, where $b= \Gamma(\alpha)\sin(\alpha \pi )/\pi$; see Figure~\ref{fig:0}.
We obtain the density function $f_\alpha(x)$ as follows
$$f_\alpha(x)=-\frac{d}{d x}E_{\alpha}(-x^\alpha)=x^{\alpha-1}E_{\alpha,\alpha}(-x^\alpha),~ 0< \alpha \leq 1,~x>0,$$
and $f_\alpha(x)=0$ for $x<0$. For $0<\alpha<1$ the asymptotic relations can be taken from \cite{Erdelyi: Magnus:Oberhettinger:Tricomi:1955} as
$$E_{\alpha}(-x^\alpha)\sim\frac{\sin(\alpha\pi)}{\pi}\frac{\Gamma(\alpha)}{x^{\alpha}},~f_\alpha(x)\sim
\frac{\sin(\alpha\pi)}{\pi}\frac{\Gamma(\alpha+1)}{x^{\alpha+1}}~ \text{for}~ x\rightarrow\infty.$$
The second asymptotic relation also comes out by formal differentiation of the first. For $0<x\rightarrow0$ we have
 $F_{\alpha}(x)=({x^{\alpha}}/{\Gamma(\alpha+1)})+\text{smaller~ order~ terms},\\
 f_\alpha(x)\sim{ x^{\alpha-1}}/{\Gamma(\alpha)}.$

 Already in the sixties of the past century \cite{Gnedenko and Kovalenko:1968} discovered our Mittag-Leffler waiting time density, $f_{\alpha}(x)$ by finding the Laplace transform of the waiting time density of a properly scaled rarefaction (thinning) limit of a renewal process with power law waiting time. But they did not identify this transform as belonging to $f_{\alpha}(x)$. In 1985 Balakrishnan found the same Laplace transform also without identifying its inverse as relevant for the time fractional diffusion process. In 1995 Hilfer and Anton were the first authors who introduced explicitly the Mittag-Leffler density $f_{\alpha}(x)$ into the theory of continuous time random walk. They showed that it is required for obtaining as evolution equation
  the fractional variant of the  Kolmogorov-Feller equation. By completely different reasoning \cite{{Mainardi:Raberto:Gorenflo:Scalas:2000}} also discussed the relevance of $f_{\alpha}(x)$ in theory of continuous time random walk. However all these early authors did not consider the renewal process with waiting time density  $f_{\alpha}(x)$ as a subject of study in its own hight but only as useful for general analysis of certain stochastic processes. The detailed investigation of the renewal process with $f_{\alpha}(x)$ as waiting time density and its analytic and probabilistic properties started (as far as we know) in 2000 with the paper by \cite{Repin and Saichev:2000}. Then more and more researchers, often independent of each other, investigated this renewal process, its properties and its applications to other process. Let us here cite  only a few relevant papers: \cite{Laskin:2003}, \cite{Mainardi:Gorenflo:Scalas:2004}, \cite{Beghin and Orsingher:2009},  \cite{Gorenflo and Mainardi:2009}, \cite{Meerschaert:Nane:Vellaisamy:2011}, \cite{Gorenflo and Mainardi:2012}.
\section{A generalized Laplacian model associated with Mittag-Leffler distribution}

\subsection{Type-2 generalized Laplacian model}
In input-output modeling, the basic idea is to model {$u=x_1-x_2$} by
 imposing assumptions on the behaviors or types of $x_1$ and $x_2$ and assumptions about whether $x_1$ and $ x_2$
 are independently varying or not, where $x_1$ and $ x_2$ respectively denote the input  and  output variables.
In a study on modeling growth-decay
mechanism, \cite{Mathai:1993} introduced a generalized Laplacian density of which the Laplace density is a special case. This concept is
connected to bilinear forms, quadratic forms and the concept of chi-squaredness of quadratic
forms, which is the basis for making inference in analysis of variance, analysis of covariance,
regression and general model building areas (see \cite{Mathai:1993}, \cite{Mathai:1993a}).

So here we introduce a contrasting growth-decay mechanism by assuming that stress and strength are
independently distributed Mittag-leffler random variables.
 Consider the random variable, $u=x_1-x_2$ which will lead to another class of generalized
 Laplacian model, say type-2 generalized Laplacian.
 The characteristic function of a type-2 generalized Laplacian model, denoted by $\phi_u(t),$
 can be obtained as follows. Here $t$ does not denote the time-variable but the argument of the
characteristic function. For $t\ge0$ we have
 \begin{eqnarray}\label{eq:3.5}
\phi_u(t)=\phi_{x_1}(t)\phi_{x_2}(-t)
&=&\frac{1}{[1+(-it)^\alpha]}\frac{1}{[1+(it)^\alpha]}\nonumber\\
&=&[1+((i)^{\alpha}+(-i)^{\alpha})t^{\alpha}+t^{2\alpha}]^{-1}\nonumber\\
&=&[1+(\rm{e}^{i(\frac{\pi \alpha}{2})}+\rm{e}^{-i(\frac{\pi \alpha}{2})})t^{\alpha}+t^{2\alpha}]^{-1}\nonumber\\
&=&\left(1+2 \cos\left(\frac{\pi \alpha}{2}\right)t^\alpha+t^{2\alpha}\right)^{-1},~0<\alpha\leq1.
\end{eqnarray}
Obviously $\phi_u(t)$ is an even function so that we finally get
\begin{equation}\label{eq:a1}\phi_u(t)=\left(1+2 \cos\left(\frac{\pi \alpha}{2}\right)|t|^\alpha+|t|^{2\alpha}\right)^{-1},t\in \mathcal{R}    ,0<\alpha\leq1.\end{equation}
The importance of this model is that we can easily obtain the fractional order residual effect.
When $\alpha=1,$ the characteristic function reduces to $(1+t^2)^{-1}$, which is the characteristic function
of a Laplace random variable whose density is $(1/2) {\exp }(-|x|)$.
 If there are several independent input variables $x_1,\ldots, x_n$ such as the situation in reaction or production problems, and if there are several independent output variables $x_{m+1},\ldots, x_{m+n}$ and if they are all independently distributed Mittag-Leffler type variables with different parameters, then the residual $v=x_1+\cdots+x_m-x_{m+1}-\cdots-x_{m+n}$ has characteristic function $$\phi_v \left(t\right)=\frac{1}{\left(1+(-it)^{\alpha_1}\right)} \cdots \frac{1}{\left(1+(-it)^{\alpha_m}\right)}\frac{1}{\left(1+(it)^{\alpha_{m+1}}\right)}\cdots \frac{1} {\left(1+(it)^{\alpha_{m+n}}\right)}, t\geq0.$$
The difference $u=x_1-x_2$ or $v=x_1+\cdots+x_m-x_{m+1}-\cdots-x_{m+n}$ can be used to describe the behaviour of a stress-strength model.
In the context of reliability the stress-strength
model describes the life of a component which has a random strength $x_1$ and is
subjected to random stress $x_2$. The component fails at the instant that the stress applied
to it exceeds the strength and the component will function satisfactorily whenever $x_1 > x_2$.
Thus $R = \text{Pr}(u > 0) = \text{Pr}(x_1 > x_2)$ is a measure of component reliability.

\subsection{Properties of type-2  generalized Laplacian distribution}
The type-2 generalized Laplacian density function, denoted by $h(u)$, can be obtained via the inverse Fourier  transform of $\phi_u(t)$. Hence
$$h(u)=\frac{1}{2\pi}\int_{-\infty}^{\infty} \frac{{\rm e}^{-itu}}{\left(1+2 \cos\left(\frac{\pi \alpha}{2}\right)|t|^\alpha+|t|^{2\alpha}\right)}{\rm d}t,0<\alpha<1, u>0.$$
\begin{prop}\label{pr2.1} For any $0<\alpha<1,$ the density function of a type-2 generalized Laplacian random variable $u$ has the representation
$$h(u)=\frac{1}{\pi}\int_0^{\infty} \frac{\cos tu}{\left(1+2 \cos\left(\frac{\pi \alpha}{2}\right)|t|^\alpha+|t|^{2\alpha}\right)}{\rm d}t, u>0.$$
\end{prop}

\subsection{Asymptotic behavior}
The anti-auto-convolution of a function vanishing for $x<0$. Assume $f(x)\equiv0$ for $x<0$. Set $g(x):=f(-x)$. Then for $h=f*g$ we get
\begin{eqnarray*}h(x)&=&(f*g)(x)\nonumber\\
&=&\int_{-\infty}^{\infty}f(y)g(x-y){\rm d} y=\int_{-\infty}^{\infty}f(y)f(y-x){\rm d} y\nonumber\\&=&\int_{y=x}^{\infty}f(y)f(y-x){\rm d} y=\int_{z=0}^{\infty}f(x+z)f(z){\rm d} z,\end{eqnarray*} and it can be shown that $(f*g)(-x)=(f*g)(x)$.
 With $f=f_{\alpha}$ we obtain (without inverting a Fourier transform) the integral representation
  $$h(x)=\int_{z=0}^{\infty}f_{\alpha}(x+z)f_{\alpha}(z){\rm d} z=
  \int_{x}^{\infty}f_{\alpha}(y)f_{\alpha}(y-x){\rm d} y,$$
  from which we can draw asymptotic relations. Because of symmetry we have $h(x)=h(-x)$ and we need only to consider $x\geq0$. For $x\rightarrow\infty$ we have with $c_{\alpha}=\{{\Gamma(\alpha+1)\sin(\alpha\pi)}\}/{\pi}$
  asymptotically $h(x)\sim\int_{y=x}^{\infty}c_{\alpha}y^{-(\alpha+1)}f_{\alpha}(y-x){\rm d}y$ and using $f_{\alpha}(x)=-{d}/{dx}( E_{\alpha}(-x^{\alpha}))$ we get by product integration
  \begin{eqnarray*}h(x)  &\sim& c_{\alpha}y^{-(\alpha+1)}E_{\alpha}(-(y-x)^{-\alpha})|_{y=x}^{y=\infty}+
  \int_{x}^{\infty}c_{\alpha}(\alpha+1)y^{-(\alpha+2)}E_{\alpha}(-(y-x)^{-\alpha}){\rm d} y\nonumber\\
  &=&c_{\alpha}x^{-(\alpha+1)}+b_{\alpha}(x).
  \end{eqnarray*}
  Because $0<E_{\alpha}(-(y-x)^{-\alpha})\leq E_{\alpha}(0)=1$ we conclude on $b_{\alpha}(x) \leq c_{\alpha}x^{-(\alpha+1)}$ so that finally $$h(x)=O(x^{-(\alpha+1)})~\text{for}~x\rightarrow\infty.$$
\begin{rem}By Tauberian theory of asymptotics for Fourier transforms (see \cite{Feller:1971},  \cite{Gorenflo and Mainardi:2009}) we find for the tail $T(x)=\int_x^{\infty}h(z){\rm d} z$
the asymptotics  $T(x)\sim\{({\Gamma(\alpha)\sin(\alpha\pi)})/{\pi}\}x^{-\alpha}$ from which by formal differentiation  we would get
 $h(x)=c_{\alpha}x^{-(\alpha+1)}$. However, aymptotic relations generally can be integrated, but differentiation needs additional smoothness requirements.
 Also we can show that $b_{\alpha}(x)$ tends to zero faster than $x^{-(\alpha+1)}$ because the estimate $0<E_{\alpha}(-(y-x)^{-\alpha})\leq E_{\alpha}(0)=1$ is very rough and the Mittag-Leffler expression tends to zero. Anyway, we now know that $h(x)=O(x^{-(\alpha+1)})~\text{for}~x\rightarrow\infty$.\end{rem}

\subsection{Moments of type-2 Laplacian distribution}
\noindent The moment $M_n=\int_{-\infty}^{\infty} x^n h(x){\rm d}x$
are given via the values at $0$ of the $n^{th}$ derivative of the Fourier transform
$${\left(1+2 \cos\left(\frac{\pi \alpha}{2}\right)|t|^\alpha+|t|^{2\alpha}\right)}^{-1}=
1-2 \cos\left(\frac{\pi \alpha}{2}\right)|t|^\alpha+\text{smaller~ order~ terms}.$$
Because $\alpha<1$ this Fourier transform is not differentiable at $t=0$. Hence, the moment $M_n$ does not exist for
$n\geq 1$. Clearly, the median exists and because of symmetry is at $0$. Contrastingly, in the limiting case $\alpha=1$ all moments exist. Then we have
$h(x)=({1}/{2}){\exp}(-|x|)$ and $M_{\beta}=\Gamma(\beta+1)$ for all real $\beta=0$.

\section{Time series model associated with the type-2 generalized Laplacian model}

\subsection{First order autoregressive model associated with the type-2 generalized Laplacian model}
Gaver and Lewis derived the exponential solution of the first order autoregressive (abbreviated as AR(1))
equation $x_n= \rho x_{n-1}+\epsilon_n, n=0,\pm1,\pm2,\cdots ,$ where
$\{\epsilon_n\}$ is a sequence of independently and identically
distributed random variables when $0\leq \rho<1$, see \cite{Gaver and Lewis:1980}.
\begin{defn}\label{df3.1} A characteristic function $\phi$ is self decomposable (belongs to class $\mathcal{L}$) if, for every $\rho, ~0<\rho<1,$ there exists a characteristic function $\phi_{\rho}$ such that $\phi(t)=\phi(\rho t)\phi_{\rho}(t),\forall t\in R$.
\end{defn}

\begin{thm}\label{th3.1} The type-2 generalized Laplacian distribution belongs to class $\mathcal{L}$.
\end{thm}

{\bf Proof.}
The proof is obvious from (\ref{eq:3.8})-(\ref{eq:3.10}).

In \cite{Gaver and Lewis:1980} it is proved that only class $\mathcal{L}$ distributions can be marginal distributions of a first order auto regressive process. Hence from Theorem \ref{th3.1} it follows that the type-2 generalized Laplacian distribution can be the marginal distribution of an AR(1) process.

The type-2 generalized Laplacian first order autoregressive process is constituted by $\{u_n;n\geq1\}$ where the $u_n$ with some $0< \rho \leq 1$
satisfy the equation
{\begin{equation}\label{eq:3.6}u_n=\rho
u_{n-1}+\epsilon_n,\end{equation}}and $\{\epsilon_n\}$ is sequence of independently and identically distributed random variables
 such that $u_n$ is stationary Markovian with type-2 generalized Laplacian distribution.
In terms of characteristic function, (\ref{eq:3.6}) can be given as
\begin{equation}\label{eq:3.7}\phi_{u_{n}}(t)=\phi_{\epsilon_{n}}(t)\phi_{u_{n-1}}(\rho t).\end{equation}
Assuming stationarity we have,
\begin{eqnarray}\label{eq:3.8}
\phi_{\epsilon_{n}(t)}&=&
\frac{\phi_{u}(t)}{\phi_{u}(\rho t)}\\
&=&\frac{(1+(-i\rho t)^{\alpha})}{(1+(-i t)^{\alpha})}\frac{(1+(i\rho t)^{\alpha})}{(1+(-i t)^{\alpha})}\nonumber\\
&=&\left[\rho^{\alpha}+(1-\rho^{\alpha})\frac{1}{(1+(-it)^{\alpha})}\right]
\left[\rho^{\alpha}+(1-\rho^{\alpha})\frac{1}{(1+(it)^{\alpha})}\right].
\end{eqnarray}
The distribution of innovation sequence can be obtained as \begin{equation}\label{eq:3.10}\epsilon_n\stackrel{d}= ML_1-ML_2,\end{equation} where
\begin{eqnarray*}
ML_1=\left\{\begin{array}{ll}
 0,& \text{with~ probability}~\rho^{\alpha}\\
 ML_{11},&\text{with~ probability}~(1-\rho^{\alpha})
 \end{array} \right.
\end{eqnarray*}
\begin{eqnarray*}
ML_2=\left\{\begin{array}{ll}
 0,& \text{with~ probability}~\rho^{\alpha}\\
 ML_{21},&\text{with~ probability}~(1-\rho^{\alpha}),
 \end{array} \right.
\end{eqnarray*}
and $ML_{11}$ and $ML_{21}$ are independently distributed Mittag-Leffler random variables.
\begin{rem}\label{re3.1}
If $u_0\stackrel{d}= ML_1-ML_2,$ then the process is strictly stationary.
\end{rem}
{\bf Proof.} For the process to be strictly stationary, it suffices to verify that
$u_n\stackrel{d}= ML_1-ML_2$  for every n. This can be proved using an inductive argument.
Suppose $u_{n-1}\stackrel{d}= ML_1-ML_2,$ then from (\ref{eq:3.5}), (\ref{eq:3.7}) and (5),
$$\phi_{u_n}(t)=\left(1+2 \cos\left(\frac{\pi \alpha}{2}\right)|t|^\alpha+|t|^{2\alpha}\right)^{-1},
  0<\alpha\leq1.$$
Hence the process is strictly stationary and Markovian, provided $u_0$ is distributed as type-2  generalized Laplacian.
\begin{rem}\label{re3.2}If $u_0$ is distributed arbitrarily and $0<\rho<1$, then the process is also asymptotically Markovian with type-2 generalized Laplacian distribution, provided $\epsilon$ is as in (\ref{eq:3.10}).\end{rem}
{\bf Proof.}
$u_n=\rho^n u_0+\sum_{k=0}^{n-1}\rho^k \epsilon_{n-k}.$
In terms of characteristic function it can be rewritten as,
$$\phi_{u_n}(t)=\phi_{u_0}(\rho^n t)\displaystyle{\prod_{k=0}^{n-1}}\phi_{\epsilon}(\rho^k t).$$
Thus the left hand side tends to
$\left(1+2 \cos\left(\frac{\pi \alpha}{2}\right)|t|^\alpha+|t|^{2\alpha}\right)^{-1},
 ~0<\alpha\leq1,$ as $n$ tends to $\infty$. Hence it follows that, even if $u_0$ is arbitrarily distributed , the process is asymptotically stationary Markovian with type-2 generalized Laplacian marginals. Thus the following theorem holds.

\begin{thm}\label{th3.2} The AR(1) process $u_n=\rho
u_{n-1}+\epsilon_n,~\rho \in(0,1)$ is strictly stationary with type-2 generalized Laplacian marginal distributions, if and only if $\{\epsilon_n\}$ are independently and identically distributed as defined in (\ref{eq:3.10}) provided $u_0$ follows a type-2 generalized Laplacian and is independent of $\epsilon_1$.\end{thm}
\subsection{Distribution of sums and joint distribution of $(u_n, u_{n+1})$}\label{subsection3}
When a stationary sequence $u_n$ is used, the distribution of the gliding sums $s_r = u_n + u_{n+1} + \cdots +
u_{n+r-1}$ is important. We have
$$u_{n+j}=\rho^j u_n+\rho^{j-1}\epsilon_{n+1}+\rho^{j-2}\epsilon_{n+2}+\cdots+\epsilon_{n+j}.$$
Hence
\begin{eqnarray*}
s_r&=&
u_n + u_{n+1} + \cdots +
u_{n+r-1}\nonumber\\
&=&\sum_{j=0}^{r-1}[\rho^j u_n+\rho^{j-1}\epsilon_{n+1}+\rho^{j-2}\epsilon_{n+2}+\cdots+\epsilon_{n+j}]\nonumber\\
&=&u_n\left( \frac{1-\rho^r}{1-\rho}\right)+\sum_{j=0}^{r-1}\epsilon_{n+j}\left( \frac{1-\rho^{r-j}}{1-\rho}\right).
\end{eqnarray*}
The characteristic function of $s_r$ is
$$\phi_{s_r}(t)= \phi_{u_n}\left( t\frac{1-\rho^r}{1-\rho}\right)\prod_{j=1}^{r-1}\phi_{\epsilon}\left( t\frac{1-\rho^{r-j}}{1-\rho}\right);$$
\begin{eqnarray*}
\phi_{s_r}(t)&=&K\left[\rho^{\alpha}+(1-\rho^{\alpha})\frac{1}{\left( 1+\left( -it\frac{1-\rho^{r-j}}{1-\rho}\right)^{\alpha}\right)}\right]
\left[\rho^{\alpha}+(1-\rho^{\alpha})\frac{1}{\left( 1+\left( it\frac{1-\rho^{r-j}}{1-\rho}\right)^{\alpha}\right)}\right],
\end{eqnarray*}
where $$K=\left[1+\left( -it\frac{1-\rho^r}{1-\rho}\right)^{\alpha}\right]^{-1}\left[1+\left( it\frac{1-\rho^r}{1-\rho}\right)^{\alpha}\right]^{-1}.$$
The density function of $s_r$ can be obtained by inverting the characteristic function as the characteristic
function uniquely determines the distribution of a random variable. Now the joint
distribution of $(u_n, u_{n+1})$ can be given in terms of characteristic function as
\begin{eqnarray}\label{eq:3.11}
\phi_{u_n,u_{n+1}}(t_1,t_2)&=&E[\rm{e}^{(it_1u_n+it_2 u_{n+1})}]\nonumber\\
&=&E[\rm{e}^{(it_1u_n+it_2 (\rho u_n+\epsilon_{n+1}))}]\nonumber\\
&=&E[\rm{e}^{(i(t_1+\rho t_2)u_n+it_2 \epsilon_{n+1})}]\nonumber\\
&=&\phi_{u_n}(t_1+\rho t_2)\phi_{\epsilon_{n+1}}(t_2)\nonumber\\
&=&I\left[\rho^{\alpha}+(1-\rho^{\alpha})\frac{1}{(1+(-it_2)^{\alpha})}\right]
\left[\rho^{\alpha}+(1-\rho^{\alpha})\frac{1}{(1+(it_2)^{\alpha})}\right],
\end{eqnarray}
where $$I={\left(1+2\cos(\frac{\pi\alpha}{2})(t_1+\rho t_2)^{\alpha}+(t_1+\rho t_2)^{2\alpha}\right)}^{-1}, t_1,t_2\geq0.$$
The above characteristic function is not symmetric in $t_1$ and $t_2$ and hence the process is not
time reversible.
\subsection{ Regression behaviour of type-2 generalized Laplacian process}\label{subsection4}
Now we shall consider the regression behaviour of the type-2 generalized Laplacian model.
Study of the regression of the model is in effect for forecasting of the process. Regression in
the forward direction explains the forecasting of future values $u_n$ while the prediction of past
values of $u_n$ can be done through regression in the backward direction. As stated in \cite{Lawrance:1978} the practical implication of regression will be in the statistical analysis of direction dependent
data, since the type-2 generalized Laplacian process is not time reversible.
\subsubsection{Regression in forward direction. }\label{subsubsection1}
The regression in the forward direction is linear since, $E[u_n|u_{n-1} = u] = \rho u, 0<\rho<1.$ Furthermore,
the conditional variance is constant.
\subsubsection{Regression in backward direction. }\label{subsubsection2}
In the backward direction, the conditional distribution of $u_n$ given $u_{n+1} = u$ has non-linear
regression and non-constant conditional variance. Following the steps described in \cite{Lawrance:1978}, the joint characteristic function of $u_n$ and $u_{n+1}$ can be derived as,
$$\phi_{u_n,u_{n+1}}(t_1,t_2)=\frac{\phi_u(t_1+\rho t_2)\phi_u(t_2)}{\phi_u(\rho t_2)}.$$
Differentiating this with respect to $t_1$ and setting $t_1 = 0, t_2 = t;$
\begin{equation}\label{eq:3.12}i E[\rm{e}^{itu_{n+1}}E[u_n|u_{n+1}]]=\frac{\phi _u^{'}(\rho t) \phi_u(t)}{\phi_u(\rho t)}=\phi _u^{'}(\rho t) \phi_{\epsilon}(t)\end{equation}
where $\phi_{\epsilon}(t)$ is as defined in (\ref{eq:3.8}), (\ref{eq:3.12}) reduces to,
\begin{equation}\label{eq:3.13}
i E[\rm{e}^{itu_{n+1}}E[u_n|u_{n+1}]]= Q\left[R +S \right],
\end{equation}
where $$Q= \frac{1}{[1+(-it)^{\alpha}][[1+(it)^{\alpha}][[1+(-i\rho t)^{\alpha}][[1+(i\rho t)^{\alpha}]},$$
$$R=\frac{-i\alpha(i \rho t)^{\alpha-1}}{[1+(- i\rho t)^\alpha][1+(i\rho t)^\alpha]^2},
S=\frac{i \alpha (-i\rho t )^{\alpha-1}}{[1+(i\rho t)^{\alpha}][1+(-i\rho t)^{\alpha}]^2}.$$
From (\ref{eq:3.13}), we can obtain the expression for $E[u_n|u_{n+1}]$ as in \cite{Lawrance:1978}. Also proceeding
with the bivariate characteristic function defined in (\ref{eq:3.11}), the conditional expectation
$E[u_n|u_{n+1}]$ can be obtained by following \cite{Gaver and Lewis:1980}.

\subsection{Simulation studies}
\subsubsection{Algorithm for  $ML_{\alpha}$ generator}
The following algorithm can be
used to generate $ML_{\alpha}$ random variables, for more details see \cite{Jayakumar and Suresh (2003)}.
\begin{enumerate}
\item Generate random variate $z$ from standard
exponential
\item Generate uniform [0,1] variate $u$, independent of $z$
\item Set $\alpha$
\item Set $w\leftarrow\sin( \pi \alpha)\cot(\pi \alpha u)-\cot(\pi \alpha )$
\item Set $y\leftarrow zw^{\frac{1}{\alpha}}$
\item Return $y$.
\end{enumerate}
We generated type-2 generalized Laplacian random variables for fixed
$\alpha=0.9$ and the histogram for those generated values
are given below.
\begin{figure}[h!]\begin{center}
\resizebox*{6.5cm}{!}{ \includegraphics{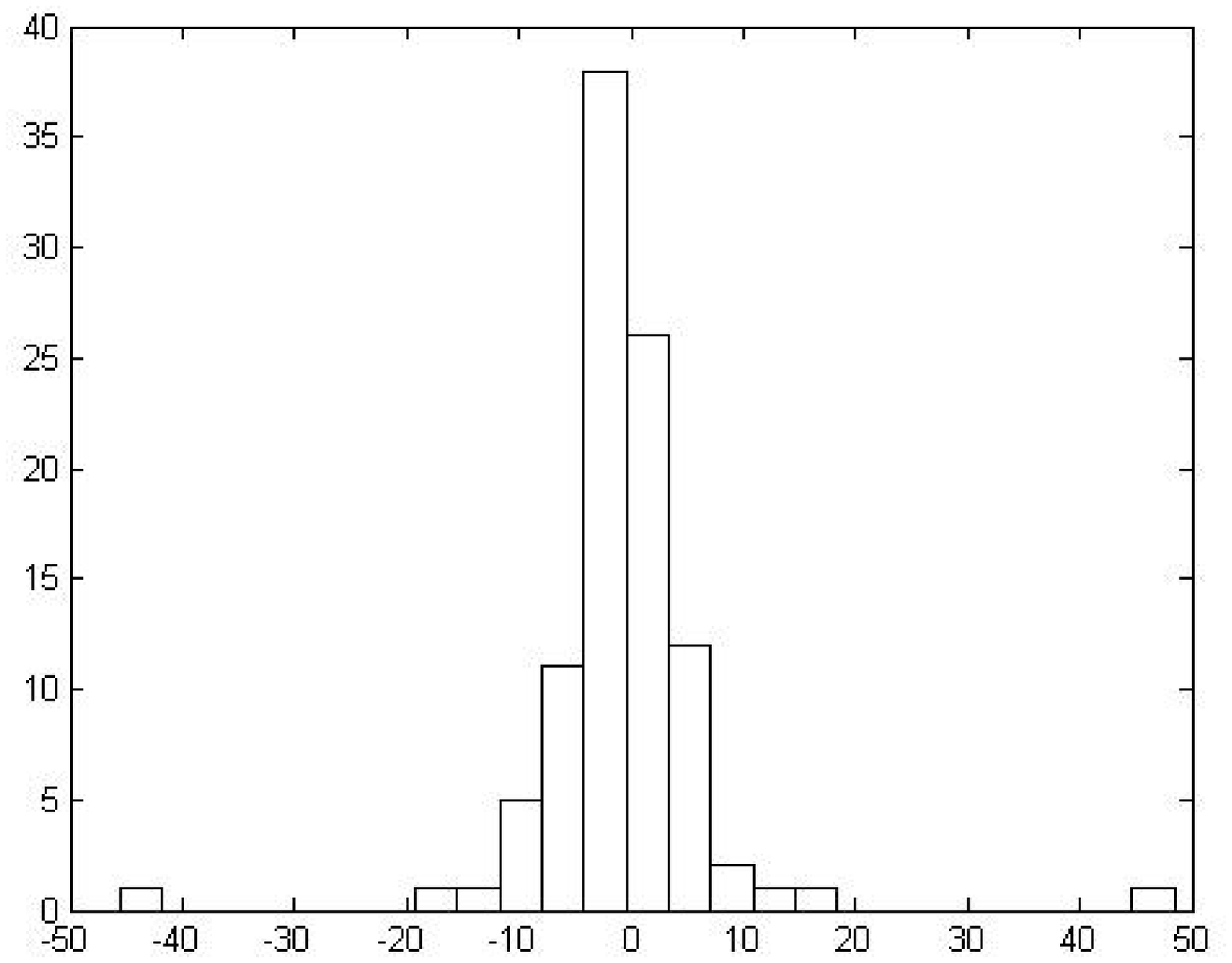}}
\end{center}\end{figure}
\subsubsection{Sample path properties }
\begin{figure}
\includegraphics[width=0.50\textwidth]{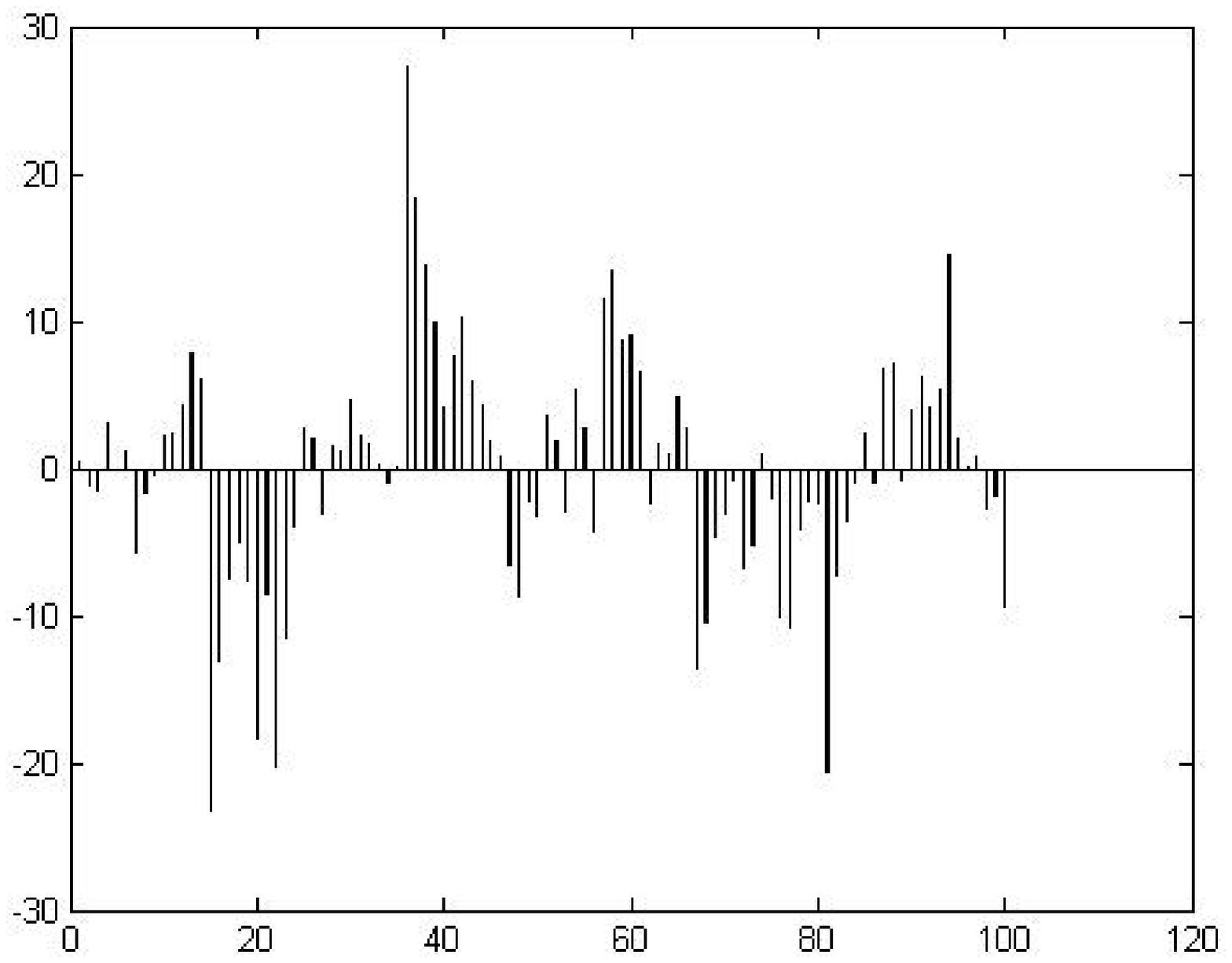} \includegraphics[width=0.50\textwidth]{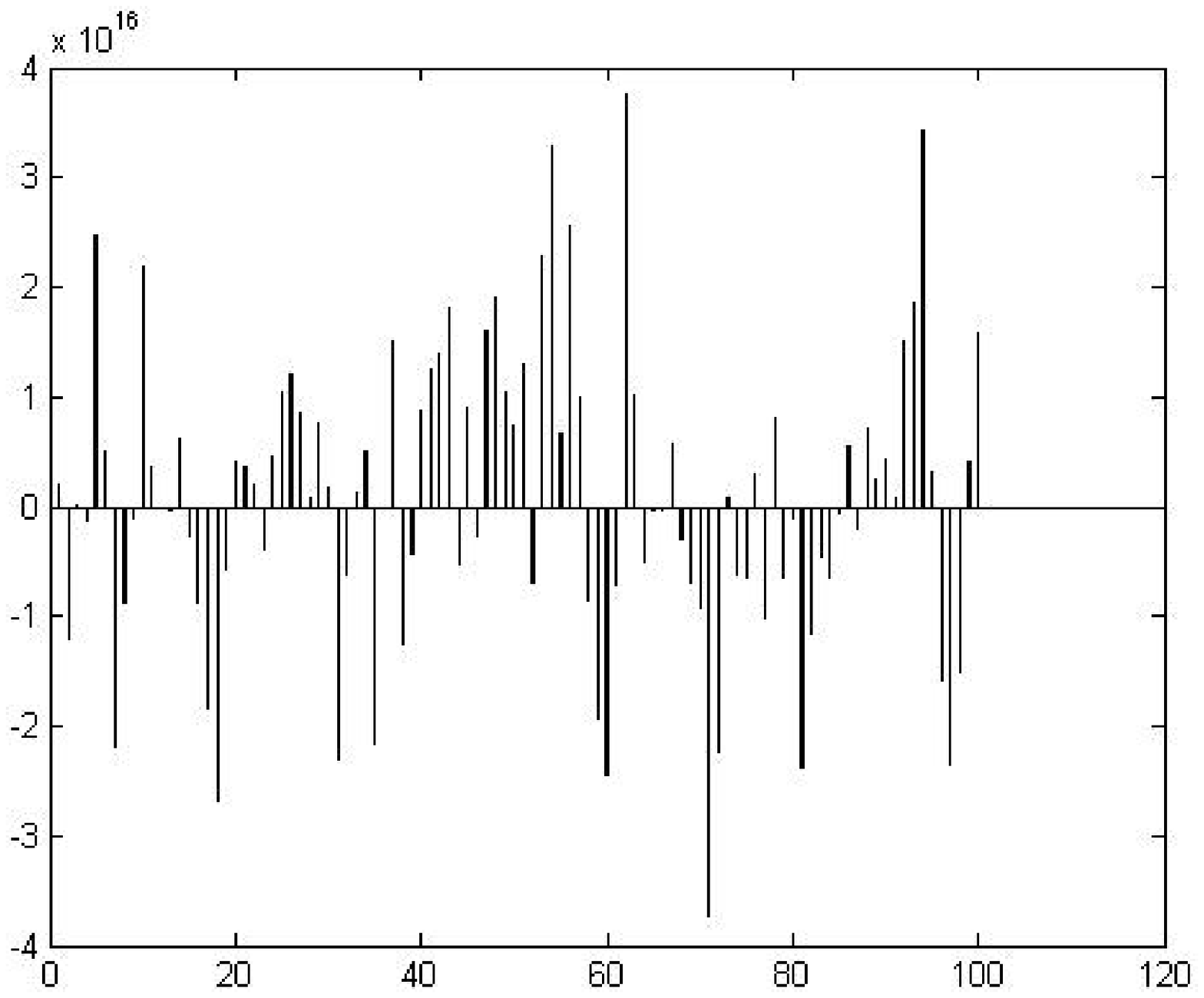}
\caption{Sample paths of type 2 generalized Laplacian process for  $\rho=0.6$ and $\alpha=0.6,0.9 $.}
\label{fig:2}
\end{figure}

Here we use the generated type-2 generalized Laplacian distribution for different values of the parameters. Its sample path is observed in the following figures. In Figure~\ref{fig:1}, we fixe $\rho= 0.3$ and the $\alpha$ values are 0.3 and 1 respectively.  For $\rho= 0.6$, we choose the $\alpha$ values as 0.6 and 0.9 respectively, the plot is given in Figure~\ref{fig:2}. It is evident from the figures that the process exhibits both positive and negative values with upward as well as downward trend. These figures point out the rich variety of contexts where the newly developed time series models can be applied. It is clear that the model gives rise to a wide variety of sample paths so that it can be used to model data from various contexts such as communication engineering, growth-decay mechanism, crop prices etc.
\begin{figure}[h!]
\includegraphics[width=0.50\textwidth]{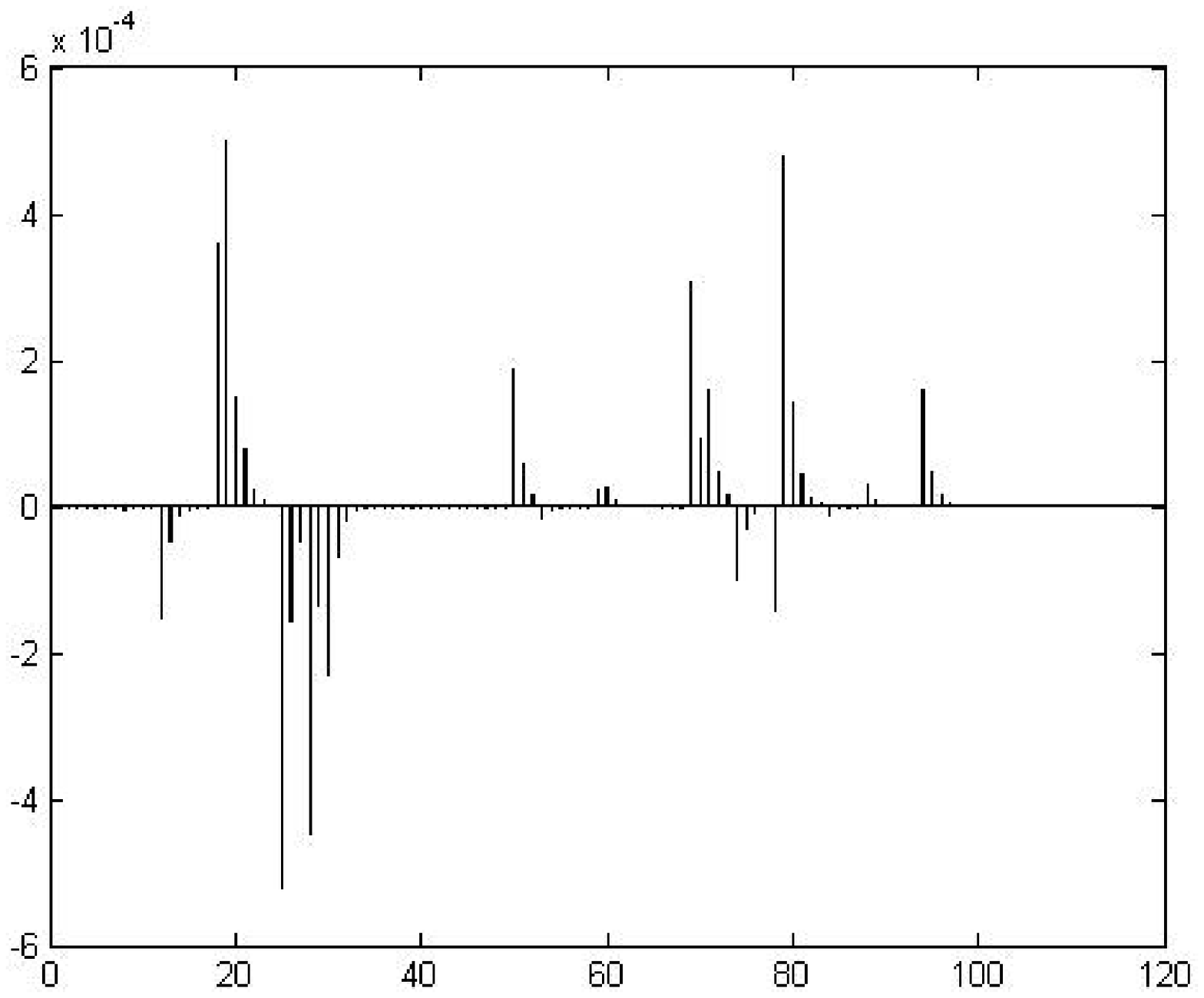} \includegraphics[width=0.50\textwidth]{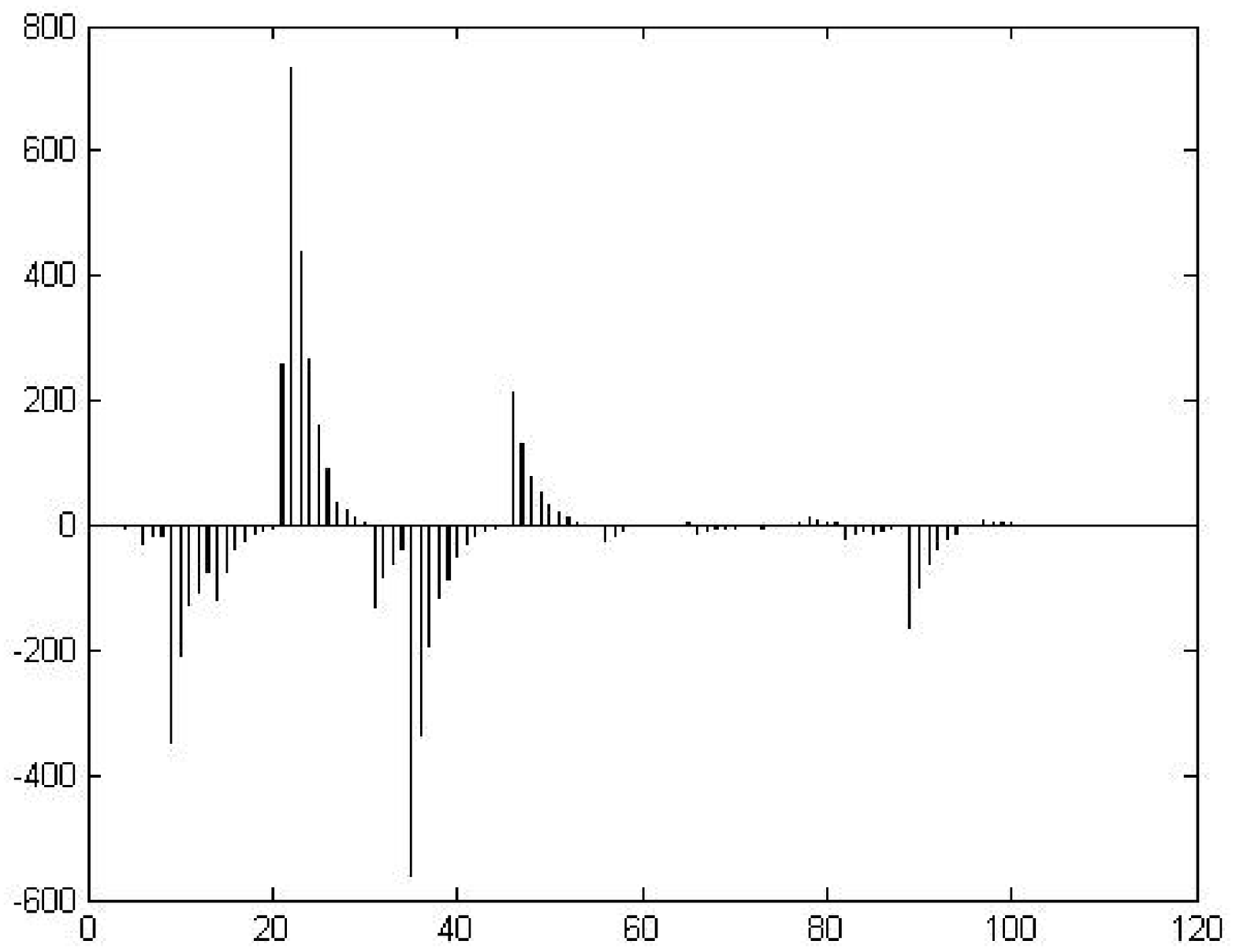}
\caption{Sample paths of type 2 generalized Laplacian process for  $\rho=0.3$ and $\alpha=0.3,1 $.}
\label{fig:1}
\end{figure}

\section{$q$-Mittag-Leffler distribution }\label{section4}
 Recently various authors have introduced several $q$-type distributions such as  $q$-exponential, $q$-Weibull, $q$-logistic and various pathway models in the context of information theory, statistical mechanics, reliability modeling etc. The $q$-exponential distribution can be viewed as a stretched model for exponential distribution so that the exponential form can be reached as $q \rightarrow 1.$ The $q$-exponential distribution is characterized by the density function
$$
f(y)=c[1+(q-1)\lambda y)]^{-\frac{1}{q-1}}; \lambda>0, y>0
$$
where $c$
is the normalizing constant. In 2010 Mathai considered the Mittag-Leffler density, associated with a
 Mittag-Leffler function as follows \cite{Mathai:2010}:
\begin{eqnarray}\label{eq:3.17}
f(y)&=&\sum_{k=0}^{\infty}\frac{(-1)^k (\eta)_k}{k!\Gamma(\alpha\eta+\alpha k)}\frac{y^{\alpha\eta-1+\alpha k}}
{(a^{\frac{1}{\alpha}})^{\alpha\eta+\alpha k}},\eta>0,a>0, 0<\alpha\leq1,0\leq y< \infty\nonumber\\
&=&\frac{y^{\alpha\eta-1}}{a^\eta}E_{\alpha,\alpha \eta}^\eta (-\frac{y^\alpha}{a}).
\end{eqnarray}
 The Laplace transform of $f$ is,
\begin{eqnarray*}
L_f(t)&=&\sum_{k=0}^{\infty}\frac{(-1)^k (\eta)_k}{k!a^{\eta+k}}\int_0^\infty\frac{y^{\alpha\eta+\alpha k-1}{\rm e}^{-ty}}{\Gamma(\alpha\eta+\alpha k)}{\rm d}y \nonumber\\
&=&\sum_{k=0}^{\infty}\frac{(-1)^k (\eta)_k}{k!a^{\eta+k}}t^{-\alpha\eta-\alpha k}=[1+at^\alpha]^{-\eta},|at^\alpha|<1.
\end{eqnarray*}
If $\eta$ is replaced by ${\eta}/{(q-1)}$ and $a$ by $a(q-1)$ with $q>1$ then we have a Laplce transform
\begin{equation}\label{eq:3.14}L_g(t)=
[1+a(q-1)t^\alpha]^{-\frac{\eta}{q-1}},q>1, t>0.\end{equation}
The distribution with Laplace transform (\ref{eq:3.14}) will be called $q$-Mittag-leffler distribution and is denoted by $ML(\alpha,\eta, q-1)$.
If $q\rightarrow 1_{+}$ in (\ref{eq:3.14}), then
$$
 \lim_{q\rightarrow 1_{+}} L_g(t)=\lim_{q\rightarrow
 1_{+}}[1+a(q-1)t^{\alpha}]^{-{{\eta}\over{q-1}}}={\rm
 e}^{-a\eta~ t^{\alpha}}=L_{f_1}(t).
$$
which is the Laplace transform of a constant multiple of a positive L\'{e}vy variable with parameter $\alpha,0<\alpha\le 1, t\geq0$. Thus  here $q$ creates a pathway of going from the general Mittag-Leffler density $f$ to a positive L\'{e}vy density $f_1$ with parameter $\alpha,$ the multiplying constant being $(a\eta)^{{1}/{\alpha}}.$

\section{$q$-Mittag-Leffler process }\label{subsection6}
\vskip.2cm
The $q$-Mittag-Leffler first order autoregressive process is constituted by $\{y_n;n\geq1\}$ where $y_n$
satisfies the equation
{\begin{equation}\label{eq:3.15}y_n=\rho
y_{n-1}+\epsilon_n,~0< \rho \leq 1,\end{equation}}where $\{\epsilon_n\}$ is sequence of independently and identically distributed random variables
 such that $y_n$ is stationary Markovian with $q$-Mittag-Leffler distribution.
We consider the $AR(1)$ structure given by (\ref{eq:3.7}). In terms of Laplace transforms, this can be rewritten as
 $$\psi_{y_{n}}(t)=\psi_{\epsilon_{n}}(t)\psi_{y_{n-1}}(\rho t).$$
Assuming stationarity we have,
\begin{eqnarray}\label{eq:3.16}
\psi_{\epsilon_{n}(t)}&=&
\frac{\psi_{y}(t)}{\psi_{y}(\rho t)}=\frac{[1+(q-1) t^{\alpha}]^{-\frac{\eta}{q-1}}}{[1+(q-1) \rho^{\alpha } t^{\alpha}]^{-\frac{\eta}{q-1}}}\nonumber\\
&=&\left[\frac{1+(q-1) \rho^{\alpha } t^{\alpha}}{1+(q-1) t^{\alpha}}\right]^{\frac{q-1}{\eta}}\nonumber\\
&=&\left[\rho^{\alpha}+(1-\rho^{\alpha})\frac{1}{[1+(q-1) t^{\alpha}]}\right]^{\frac{q-1}{\eta}}, t>0.
\end{eqnarray}

The infinitely divisible $ML(\alpha,\eta, q-1)$ variable is of the class $\mathcal{L}$, and therefore that (\ref{eq:3.16}) is the Laplace-Stieltjes transform of a distribution function follows from the class $\mathcal{L}$ theorem of \cite{Feller:1971}, since the determining canonical measure M of the $ML(\alpha,\eta, q-1)$ variable is $q-1$ (when $\eta=1$) times the determining canonical measure of the $ML(\alpha,\eta, q-1)$ variable. Thus we can in principle generate an autoregressive process with gamma marginals (\ref{eq:3.17}) by utilizing the $\{\epsilon_n\}$ process characterized by (\ref{eq:3.16}). Here are three simple special cases.
When $q=2$ and $\eta=1$ in (\ref{eq:3.16}) we have
\begin{equation}\label{eq:3.18}\psi_{\epsilon_1}(t)=\left[\rho^{\alpha}+(1-\rho^{\alpha})\left(\frac{1}{1+ t^{\alpha}}\right)\right], t>0.\end{equation}
Thus the random variable,
\begin{eqnarray*}
\epsilon_1=\left\{\begin{array}{ll}
 0,&\text{with~ probability}~\rho^{\alpha}\\
 M,&\text{with~ probability}~(1-\rho^{\alpha}).
 \end{array} \right.
\end{eqnarray*}
Hence, $\epsilon_1$ is a convolution of an atom of mass $\rho^{\alpha}$ at zero and $ 1-\rho^{\alpha}$ at $M$ where $M$ is distributed as $ML(\alpha)$.
When $q=3$ and $\eta=1$ in (\ref{eq:3.16}) we have
\begin{equation}\label{eq:3.19}\psi_{\epsilon_2}(t)=\left[\rho^{2\alpha}+2\rho^{\alpha} (1-\rho^{\alpha})\left(\frac{1}{1+ 2t^{\alpha}}\right)+ (1-\rho^{\alpha})^2\left(\frac{1}{1+ 2t^{\alpha}}\right)^2\right], t>0,\end{equation}
\begin{eqnarray*}
\text{and}~\epsilon_2=\left\{\begin{array}{ll}
 0,&\text{with~ probability}~\rho^{2\alpha}\\
 ML(\alpha,3, 2),&\text{with~ probability}~2\rho^{\alpha}(1-\rho^{\alpha})\\
 ML(\alpha,4, 2),&\text{with~ probability}~(1-\rho^{\alpha})^2.
 \end{array} \right.
\end{eqnarray*}
When $q={3}/{2}$ and $\eta={1}/{4}$ in (\ref{eq:3.16}) we have
\begin{equation}\label{eq:3.19}\psi_{\epsilon_3}(t)=\left[\rho^{2\alpha}+2\rho^{\alpha} (1-\rho^{\alpha})\left(\frac{1}{1+ \frac{t^{\alpha}}{2}}\right)+ (1-\rho^{\alpha})^2\left(\frac{1}{1+ \frac{t^{\alpha}}{2}}\right)^2\right], t>0,\end{equation}
\begin{eqnarray*}
\text{and}~\epsilon_3=\left\{\begin{array}{ll}
 0,&\text{with~ probability}~\rho^{2\alpha}\\
 ML(\alpha,{1}/{2}, {1}/{2}),&\text{with ~probability}~2\rho^{\alpha}(1-\rho^{\alpha})\\
 ML(\alpha,1, {1}/{2}),&\text{with~ probability}~(1-\rho^{\alpha})^2.
 \end{array} \right.
\end{eqnarray*}

In general the $q$-Mittag-Leffler process can give a generalization of the model given in \cite{Gaver and Lewis:1980}. Hence the essentials of fractional calculus according to different approaches that can be useful for our applications in the theory of probability and stochastic processes are established.

\section{ Applications }\label{section5}
During the last 15 years a lot of engineers and scientists have shown very much interest in the Mittag-Leffler function and
Mittag-Leffler type functions due to
their vast potential of applications in several fields such as fluid flow, rheology, electric networks, probability, and statistical distribution
theory. The Mittag-Leffler function arises naturally in the solution of fractional order integral
or differential equations, and especially in the investigations of the
fractional generalization of the kinetic equation, random walks, L\'{e}vy flights, anomalous diffusion
transport and in the study of complex systems. In recent years the fractional  generalization of the classical Poisson process has
    gained increasing interest. In it the waiting time between events is the
    Mittag-Leffler distribution function $F_\alpha$ in place of the exponential
    distribution. Of all the papers devoted to this special renewal process we
    content ourselves to cite only \cite{Mainardi:Raberto:Gorenflo:Scalas:2000}, \cite{Mainardi:Gorenflo:Scalas:2004} and \cite{Gorenflo and Mainardi:2012} . Mittag-Leffler distributions can be used as waiting-time
distributions as well as first-passage time distributions for
certain renewal processes with geometric exponential as waiting-time  distribution. They can also be used in reliability modeling as an alternative
for exponential lifetime distribution.
The ordinary and generalized Mittag-Leffler
functions interpolate between a purely exponential law and power-law-like behavior of
phenomena governed by ordinary kinetic equations and their fractional counterparts, see
\cite{Mathai:Saxena:Haubold:2006}.
Mittag-Leffler functions are also used for
computation of the change of the chemical composition in stars like
the Sun. Recent investigations have proved that they are useful in modelling  the flux of solar neutrinos in cosmological studies, which can be expressed in terms of special
functions like $G$ and $H$-functions, see \cite{Mathai:Saxena:Haubold:2006},  \cite{Sebastian:2011}. The Mittag-Leffler distribution finds applications in a wide range of contexts such as stress-strength analysis, growth-decay mechanisms like formation of
sand dunes in nature, input-output situations in economics,  industrial productions,
production of melatonin in human body etc.\\

\section{Acknowledgements}
First author concedes gratefully all mentors for their comments and constructive suggestions as it helped to improve the article.


\begin{thebibliography}{12}
%




\bibitem{Balakrishnan:1985}
{Balakrishnan, V.}, {Anomalous diffusion in one dimension}, {Physica A}, {132}, 569-580 (1985).


\bibitem{Beghin and Orsingher:2009}
{Beghin, L.} and  {Orsingher, E.},  {Fractional Poisson processes and related
random motions}, {Electronic Journ. Prob.}, {{14(61)}}, 1790-1826 (2009).


\bibitem{Erdelyi: Magnus:Oberhettinger:Tricomi:1955}
{Erd\'{e}lyi, A., Magnus, W., Oberhettinger, F.}  and  {Tricomi, F. G.}, { {Higher Transcendental Functions}}, vol. 3, McGraw-Hill, New York (1955).

\bibitem{Feller:1971}
{Feller, W.}, {{An Introduction to Probability Theory and Its Applications}}, vol.2, Wiley, New York (1971).

\bibitem{Gaver and Lewis:1980}
{Gaver, D. P.} and {Lewis, P. A. W.}, {First-order autoregressive gamma sequences and point processes}, {Adv. Appl.Prob.}, {{12}}, 727-745 (1980).



\bibitem{Gnedenko and Kovalenko:1968}
{Gnedenko , B.V.} and  {Kovalenko, I.N.}, {{Introduction to Queueing Theory,
Israel Program for Scientific Translations}}, Jerusalem (1968).

\bibitem{Gorenflo and Mainardi:1997}
{Gorenflo, R.}  and {Mainardi, F.}, {Fractional calculus: integral and differential equations of fractional order},
in A. Carpinteri and F. Mainardi (Editors), Fractals and Fractional Calculus in Continuum
Mechanics,  Springer Verlag,  223-276 (1997).


\bibitem{Gorenflo and Mainardi:2009}
{Gorenflo, R.} and {Mainardi, F.}, {Some recent advances in theory and simulation of fractional diffusion processes}, {Journal of Computational and Applied Mathematics}, {229}, 400-415 (2009).

\bibitem{Gorenflo and Mainardi:2012}
{Gorenflo, R.} and {Mainardi, F.}, {Laplace-Laplace analysis of the fractional Poisson process}, {Proceedings of Analytical Methods of Analysis and Differential Equations}, Minsk, 41-56 (2012).


\bibitem{Hilfer  and Anton:1995}
{Hilfer , R.}  and {Anton, L.}, {Fractional master equations and fractal time
random walks}, {Physical Review E}, {51}, R848-R851 (1995).



\bibitem{Jayakumar and Suresh (2003)}
Jayakumar, K. and Suresh, R. P.,   Mittag-Leffler distributions, J. Indian Soc. Prob. Statist., {7}, 51-71, (2003).


\bibitem{Laskin:2003}
{Laskin, N.}, {Fractional Poisson processes}, {Communications in Nonlinear Science
and Numerical Simulation}, {8}, 201-213 (2003).


\bibitem{Lawrance:1978}
{Lawrance, A. J.}, {Some autoregressive models for point processes}, {Proceedings of Bolyai Mathematical Society Colloquium on Point Processes and Queuing Problems} Hungary, {24}, 257-275 (1978).



\bibitem{Mainardi and Gorenflo:2000}
{Mainardi, F.} and {Gorenflo, R.}, {On Mittag-Leffler-type functions in fractional evolution processes}, {Journal of Computational and Applied Mathematics}, {118}, 283-299 (2000).

\bibitem{Mainardi:Raberto:Gorenflo:Scalas:2000}
{Mainardi, F., Raberto, M., Gorenflo, R.} and {Scalas, E.}, {Fractional calculus
and continuous-time finance II: the waiting-time distribution}, {
Physica A}, {287(3-4)}, 468-481  (2000).

\bibitem{Mainardi:Gorenflo:Scalas:2004}
{Mainardi, F., Gorenflo, R.} and {Scalas, E.},  {A fractional generalization of the Poisson process}, {Vietnam Journal of Mathematics}, {32 SI}, 53-64  (2004).

\bibitem{Mathai:1993} {Mathai, A.M.}, {On non-central generalized Laplacianness of
quadratic forms in normal variables}, {Journal of Multivariate Analysis}, {45}, 239-246  (1993).


\bibitem{Mathai:1993a}
Mathai, A. M., The residual effect of growth-decay mechanism and the
distributions of covariance structure, {The Canadian Journal of
Statistics}, {21(3)},  277-283 (1993a).


\bibitem{Mathai:Saxena:Haubold:2006}
{Mathai, A. M., Saxena, R. K.}
and {Haubold, H. J.}, {A certain class of laplace transforms with applications
to reaction and reaction-diffusion equations}, {
Astrophysics \& Space Science}, {305}, 283-288 (2006).


\bibitem{Mathai:2010} {Mathai, A. M.}, {Some properties of Mittag-Leffler functions and matrix-variate analogues: A statistical
perspective}, {Fractional Calculus and Applied Analysis}, {13}, 113-132  (2010).

\bibitem{Meerschaert:Nane:Vellaisamy:2011}
{Meerschaert, M. M., Nane, E.}
 and {Vellaisamy, P.}, {The fractional Poisson
process and the inverse stable subordinator}, {
Electronic Journ. Prob.}, {16}, 1600-1620  (2011).
\bibitem{Pillai:1990}
{Pillai, R. N.}, {On Mittag-Leffler Functions and Related Distributions}, {Ann. Inst.Statist.Math.}, {42(1)}, 157-161 (1990).

\bibitem{Podlubny:1999}
{Podlubny I.}, {{Fractional Differential Equations}},  Academic Press, San Diego (1999).

\bibitem{Repin and Saichev:2000}
{Repin, O.N.}  and {Saichev, A.I.}, {Fractional Poisson law}, {Radiophysics and
Quantum Electronics}, {43(9)}, 738-741 (2000).




\bibitem{Ross :1996}
{Ross, S.M.}, {{Stochastic Processes}}, 2-nd Edition, Wiley, New York (1996).



\bibitem{Sebastian:2011}
{Sebastian, N.}, {A generalized gamma model
associated with a Bessel function}, {Integral transforms and
Special functions}, {22(9)}, 631-645  (2011).












\end{thebibliography}
\end{document}